\documentclass[reqno,a4paper,12pt]{amsart} 

\usepackage{amsmath,amscd,amsfonts,amssymb,graphicx}
\usepackage{mathrsfs,dsfont}

\usepackage{tikz}
\usetikzlibrary{patterns}
\usepackage{float}
\usepackage{appendix}
\usepackage{caption}
\usepackage{subcaption}
\usepackage{enumerate}

\numberwithin{equation}{section}
\numberwithin{figure}{section}

\def\R{\mathbb{R}}

\def\Z{\mathbb{Z}}

\def\1{\mathds{1}}

\renewcommand\le{\leqslant}
\renewcommand\ge{\geqslant}
\renewcommand\leq{\leqslant}
\renewcommand\geq{\geqslant}

\theoremstyle{plain}
\newtheorem{theorem}{Theorem}[section]
\newtheorem{lemma}[theorem]{Lemma}

\newtheorem*{claim*}{Claim}
\newtheorem*{theorem*}{Theorem}

\theoremstyle{definition}
\newtheorem{definition}[theorem]{Definition}
\newtheorem*{definition*}{Definition}
\newtheorem*{remarks*}{Remarks}
\newtheorem*{remark*}{Remark}
\newtheorem{remark}[theorem]{Remark}

\newenvironment{enumerate-math}
{\begin{enumerate}
\addtolength{\itemsep}{5pt}
}
{\end{enumerate}}

\newenvironment{dedication}
  {%\clearpage           % we want a new page
   \thispagestyle{empty}% no header and footer
   \itshape             % the text is in italics
  }
  {\par % end the paragraph
  }

\newenvironment{enumerate-text}
{\begin{enumerate}
\addtolength{\itemsep}{5pt}
}
{\end{enumerate}}

\begin{document}
\title{Fourier frames for surface-carried measures }

\author{Alex Iosevich}
\email{iosevich@math.rochester.edu}
\address{Department of Mathematics, University of Rochester, Rochester, NY, 14627} 

\author{Chun-Kit Lai}
 \email{cklai@@sfsu.edu}
 \address{Department of Mathematics, San Francisco State University, San Francisco, CA 94132.}

 \author{Bochen Liu}
 \email{Bochen.Liu1989@gmail.com}
 \address{Department of Mathematics, The Chinese University of Hong Kong, Ma Liu Shui, Shatin, Hong Kong}

\author{Emmett Wyman}
\email{ewyman@math.northwestern.edu}
\address{Department of Mathematics, Northwestern University, Evanston, IL 60208}

\thanks{Iosevich is partially supported by Dean's Research Funds at the University of Rochester. Liu is partially supported by the grant CUHK24300915 from the Hong Kong Research Grant Council, and a direct grant of research (4053341) from the Chinese University of Hong Kong.}
\subjclass[2010]{Primary 42B05, 42A85.}
\keywords{Gaussian urvature, Fourier frames, surface measures, spectral measures, Riemannian manifolds. }

\begin{abstract}
In this paper we show that the surface measure on the boundary of a convex body of everywhere positive Gaussian curvature does not admit a Fourier frame. This answers a question proposed by Lev and provides the first example of a uniformly distributed measure supported on a set of Lebesgue measure zero that does not admit a Fourier frame. In contrast, we show that the surface measure on the boundary of a polytope always admits a Fourier frame.

We also explore orthogonal bases and frames adopted to sets under consideration. More precisely, given a compact manifold $M$ without a boundary and $D \subset M$, we ask whether $L^2(D)$ possesses an orthogonal basis of eigenfunctions. The non-abelian nature of this problem, in general, puts it outside the realm of the previously explored questions about the existence of bases of characters for subsets of locally compact abelian groups.
\end{abstract}

\maketitle

\begin{dedication} This paper is dedicated to Alexander Olevskii on the occasion of his birthday. Olevskii's mathematical depth and personal kindness serve as a major source of inspiration for us and many others in the field of mathematics. \end{dedication} 

\section{Introduction}

\subsection{Background.} Let $\mu$ be a finite positive Borel measure on ${\mathbb R}^d$. We say that $\mu$ is a {\it frame-spectral measure} if there exists a set of exponential functions $E(\Lambda) = \{e^{2\pi i \lambda\cdot x}: \lambda\in\Lambda\}$ such that $E(\Lambda)$ forms a {\it Fourier frame} for $L^2(\mu)$, in the sense that there exist constants $0<A\le B<\infty$ such that 
\begin{equation}\label{Fframe}
A \int |f(x)|^2d\mu(x) \le \sum_{\lambda\in\Lambda} \left| \int f(x)e^{-2\pi i \lambda\cdot x}d\mu(x)\right|^2 \le B \int |f(x)|^2d\mu(x) 
\end{equation}
for all $f\in L^2(\mu)$. If such $E(\Lambda)$ exists, we call $\Lambda$ a {\it frame-spectrum} for $\mu$. If only the upper bound holds in (\ref{Fframe}), we call $E(\Lambda)$ a {\it Bessel sequence} for $\mu$. If $L^2(\mu)$ admits a Fourier  orthonormal basis $E(\Lambda)$, we call $\mu$ a {\it spectral measure} and $\Lambda$ a {\it spectrum} for $\mu$.   A set $\Omega \subset {\Bbb R}^d$ is called a {\it spectral set} if the restriction of the Lebesgue measure to $\Omega$ is a spectral measure. It is clear that a Fourier orthonormal basis is also a Fourier frame, so a spectral measure is always frame-spectral, but not necessarily the other way around.

\medskip

Fourier frames were introduced by Duffin and Schaeffer \cite{DS} and have since become a fundamental tool in signal processing and data transmission. An interested reader can find the general frame theory background description in \cite{Chr,H} and some Fourier frame theory background in \cite{OU}. One of the key questions in this subject matter is the following: 

\medskip

{\bf (Qu 1):} Which measures $\mu$ are frame-spectral?

\medskip

This problem was first studied by Fuglede \cite{Fu}. His celebrated conjecture, which asserts that a spectral set and a translational tile are equivalent, was disproved in both direction \cite{T,KM06}, but remains captivating among many researchers up-to-date.  This problem was then advanced to singular measures by Jorgensen and Pedersen \cite{JP98}, who discovered that the middle-fourth Cantor measure is a spectral measure, while the the middle-third Cantor measure is not. Proposed first by Strichartz \cite{St1},  it is still an open problem today whether the middle-third Cantor measure is frame-spectral.

\medskip

 It is known that a compactly supported frame-spectral measure $\mu$ must be of pure type \cite{HLL}, which means that it is purely atomic, purely absolutely continuous with respect to Lebesgue measure, or purely singular with respect to the Lebesgue measure. If this measure is purely atomic, it is  known that it is frame-spectral if and only if the measure has only finitely many atoms. If this measure is absolutely continuous, then it is frame-spectral if and only if the density of the measure is bounded above and bounded away from zero almost everywhere on the essential support of $\mu$ \cite{Lai,DL14,NOU}.  

\medskip

The case for the singular measures is much less well understood. Dutkay and Lai \cite{DL14} proposed that a uniformity principle, formulated in terms of ``translational absolute continuity,"  should be a necessary condition. They showed that a self-similar measure with non-uniform probability weight cannot be frame-spectral. This problem was further studied in \cite{FL18,Lev}.  Denote by ${\mathcal H}^s$ the $s$-dimensional Hausdorff measure.  The uniformity principle leads us naturally to study the following generalization of Strichartz' question: 

\medskip

{\bf (Qu 2):} Let $K$ be a measurable set such that $0<{\mathcal H}^s (K)<\infty$. Does the measure $\mu = {\mathcal H}^s|_K$ admit a Fourier frame? 

\medskip

 If $s = d$,  {\bf (Qu 2)} is trivial for bounded set $K$ since we can put $K$ inside a cube and the exponential orthonormal basis on the cube immediately induces a tight frame on $K$. However, such constructions cease to exist if $K$ has Lebesgue measure zero or $K$ is unbounded. However, when $K$ is unbounded with positive finite Lebesgue measure, Nitzan, Olevskii and Ulanoskii \cite{NOU} showed that $L^2(K)$ also admits a Fourier frame with the help of the Kadison-Singer theorem \cite{MSS}. When the set has Lebesgue measure zero,  to the best of our knowledge, all measures of the form  ${\mathcal H}^s|_K$ either admits a Fourier frame or were not known if there exists any Fourier frames.

\subsection{Main Results} 
The purpose of the present paper is to demonstrate, for the first time, a measure of the form  ${\mathcal H}^s|_K$ that does not admit a Fourier frame. In the process, we shall answer a question posed by Nir Lev \cite{Lev}. 

\begin{theorem}\label{maintheorem1}
Let $K$ be a convex body on ${\mathbb R}^d$ with smooth boundary $\partial K$ having everywhere positive Gaussian curvature and let $\sigma$ be the surface  measure supported on $\partial K$. Then the measure $\sigma$ does not admit a Fourier frame.
\end{theorem}

\medskip

\begin{remark} Let $\mu_{\delta}$ denote $\delta^{-1}$ times the indicator function of the annulus of radius $1$ and thickness $\delta$ in ${\Bbb R}^d$. Since the annulus has positive Lebesgue measure, the aforementioned result due to Nitzan, Olevskii and Ulanovskii (\cite{NOU}) implies that there exist universal constants $C,c$ (independent of $\delta$) and a set $\Lambda_{\delta}$ such that for $f \in L^2(\mu_{\delta})$ 
$$ c{||f||}_{L^2(\mu_{\delta})}^2 \leq \sum_{\lambda \in \Lambda_{\delta}} {|\widehat{f d\mu_{\delta}}(\lambda)|}^2 \leq C {||f||}_{L^2(\mu_{\delta})}^2.$$

In the usual sense, $\mu_{\delta} \to \sigma$, the surface measure on the unit sphere $S^{d-1}$. Nevertheless, Theorem \ref{maintheorem1} implies that $L^2(\sigma)$ does not possess a frame of exponentials. This shows that the usual weak limit arguments that appeared in \cite{NOU} or other papers will not work well in constructing Fourier frames in the singular measure setting.
\end{remark} 

Theorem \ref{maintheorem1} will follow from the two general theorems about the summability of the Fourier frame spectra. We will prove that the surface measure on $\partial K$ will satisify both conditions stated in the theorems below with $\gamma = d-1$, which leads to a contradiction.  % For the precise definition of Beurling density and Beurling/Fourier dimension, one may refer to Section 2. 

\begin{theorem}\label{thdecay}
Let $\mu$ be a finite Borel measure on ${\mathbb R}^d$, $\Lambda\subset {\mathbb R}^d$ be a countable set such that for some $A>0$,
\begin{equation}\label{frame_lower_bound}A \int |f(x)|^2d\mu(x) \le \sum_{\lambda\in\Lambda} \left| \int f(x)e^{-2\pi i \lambda\cdot x}d\mu(x)\right|^2,\ \ \ \forall\,f\in L^2(\mu). \end{equation}
Suppose there exists $C>0$, $0<\gamma\leq d$ such that
\begin{equation}\label{decay_eq}
|\widehat{\mu}(\xi)|\le C |\xi|^{-\gamma/2},\ \ \forall\, |\xi|>1.
\end{equation}
Then
$$
\sum_{\lambda\in\Lambda\setminus\{0\}}\frac{1}{|\lambda|^{\gamma}} = \infty.
$$
\end{theorem}

\vskip.125in 

\begin{theorem}\label{thfinite}
Let $\mu$ be a finite Borel measure on ${\mathbb R}^d$, $\Lambda\subset {\mathbb R}^d$ be a countable set such that $E(\Lambda)$ is a Bessel sequence of $\mu$, namely for some $B>0$
$$\sum_{\lambda\in\Lambda} \left| \int f(x)e^{-2\pi i \lambda\cdot x}d\mu(x)\right|^2 \le B \int |f(x)|^2d\mu(x),\ \ \ \forall\,f\in L^2(\mu).$$
Suppose there exist $r, L, \gamma, c>0$ such that for any $\lambda\in\Lambda, \, |\lambda|>L$,
\begin{equation}\label{eqthfinite}
|\lambda|^{\gamma} \int_{B_r(\lambda)} |\widehat{\mu}(\xi)|^2\,d\xi > c >0.
\end{equation}
Then 
$$
\sum_{\lambda\in\Lambda\setminus\{0\}}\frac{1}{|\lambda|^{\gamma}} <\infty.
$$
\end{theorem}

%\vskip.125in 

%\begin{remark} Strichartz (\cite{St89}) proved that the estimate (\ref{eqthfinite}) holds for a class of measures that certainly includes the surface measure on the boundary of a convex body. For the sake of a self-contained presentation we reprove this fact, in the specal case under consideration, in the course of proving Theorem \ref{maintheorem1}. \end{remark} 

\medskip

 Theorem \ref{thdecay} is interesting in its own right. A classical result of Landau \cite{Lan} states that {\it if $\Lambda$ is a frame-spectrum for $L^2(\Omega)$, then the lower Beurling density of $\Lambda$ is at least the Lebesgue measure of $\Omega$}. This implies that $\Lambda$ is distributed like a lattice and therefore $\sum_{\lambda\in\Lambda\setminus\{0\}}|\lambda|^{-d} = \infty$ trivially. The Landau's result has produced a lot of important applications in frame theory (see e.g. \cite{Chr}).
 Unfortunately, such density result is longer true in the fractal setting. It was  found that  the standard middle-fourth Cantor measure admits an exponential orthonormal basis of frequency spectrum $\Lambda$ as sparse as we wanted \cite{DHL13}, which means the sum of some spectra could be finite for all $\gamma>0$. While it is well-known that the middle-fourth Cantor measure does not have any Fourier decay as in (\ref{decay_eq}),  we can view Theorem \ref{thdecay} as a natural generalization of the classical Landau density result to the singular measures.  
 
 \medskip

%\subsection{Surface-carried measures}  We now turn to our attention to the frame-spectrality of the  surface measure.  A folklore problem that people has been studying is the specific problem:

%\medskip

%{\bf (Qu 3):} Are the surface measure of the unit sphere  frame-spectral? 

%\medskip

%Lev showed that a spherical cap is frame-spectral and he explicitly asked ({\bf Qu 3}) in the end of his paper \cite{Lev}. In this paper, we show that the surface measure of the unit sphere is not spectral. More generally, we have

%\medskip

%We now sketch the proof for Theorem \ref{maintheorem1}. We know that the Fourier transform for $\sigma$ and $\chi_{K}$ shares a similar Fourier asymptotics. The later was studied by Iosevich and Rudnev \cite{IR}. They showed that a spectrum, if exists, must lie on a line. This is also the case for the boundary surface measure $\sigma$. If $d=2$, we parametrize a portion of the surface as a graph with the line as one of the axes. a spectrum must then be a tight Fourier frame on that surface. However, from a result in \cite{DL14}, this is not possible unless the curvature is zero. For $d\ge 3$, as it is well-known that the surface measure decays in the order to $d-1 \ge 2$,  using Theorem \ref{thdecay}, the  Beurling dimension of $\Lambda$ cannot be attained if it is just lying in a line. 

In contrast to the case of the positive Gaussian curvature, we also study the (flat) surface measure of polytopes that need not be convex. We show that they are all frame-spectral. 

\begin{theorem}\label{maintheorem2}
Let $K$ be a polytope on ${\mathbb R}^d$ and let $\sigma$ be the surface measure supported on $\partial K$. Then the measure $\sigma$  is frame-spectral.
\end{theorem}

 The proof of Theorem \ref{maintheorem2} is partly inspired by Lev's argument in \cite[Theorem 1.1]{Lev} where he proves that a sum of two singular frame-spectral measures is frame-spectral if they have no atoms and they are supported in two distinct orthogonal subspaces whose intersection is trivial. We note here that $\sigma$ can be written as a finite sum of $(d-1)$-Hausdorff measures. The subspaces they are supported on may or may not be orthogonal and they may intersect non-trivially. We can control the frame bound in this case because the $\sigma$ we consider here are the sum of  Lebesgue measures on lower dimensional subspaces as opposed to general frame-spectral measures in Lev's theorem. It would be nice if Theorem \ref{maintheorem2} can be generalized to a finite sum of arbitrary frame-spectral measures.  

\medskip

  A more general version of Theorem \ref{maintheorem2} will be proved in Theorem \ref{theorem_flat_surface}, in which some surface measures of different Hausdorff dimensions are allowed in the finite sum.

\medskip

\subsection{Orthogonal Eigenbases for Riemannian Manifold.} As we have seen, $L^2(\sigma)$ does not possess a frame of exponentials, where $\sigma$ is the surface measure on the sphere, demonstrating that exponential systems are not very efficient in this realm. However, $L^2(\sigma)$ possesses an orthogonal basis consisting of spherical harmonics. More generally, if $M$ is a compact manifold without a boundary, then $L^2(M)$ possesses naturally an orthogonal basis $\{e_j\}$, where $e_j$ are the eigenfunctions of the Laplace-Beltrami operator on $M$. This leads us to the following question: 

\medskip 

{\bf (Qu 3):} Let $M$ be a compact Riemannian manifold without a boundary. For which subsets $D \subset M$, does $L^2(D)$ possess an orthogonal basis of eigenfunctions of the Laplace-Beltrami operator? 

\medskip

 To set up notations, we let $(M,g)$ be a compact Riemannian manifold without boundary, and let $\Delta_g$ be the Laplace-Beltrami operator on $M$, given in local coordinates by
\[
	\Delta_g f = |g|^{-1/2} \sum_{i,j} \partial_i(g^{ij} |g|^{1/2} \partial_j f).
\]
The Laplace-Beltrami operator commutes with isometries, i.e. if $\phi$ is an isometry of $M$, then $\Delta_g (f \circ \phi) = (\Delta_g f) \circ \phi$. By the spectral theorem, there exists an orthonormal basis of eigenfunctions $e_j$ for $j = 1,2,\ldots$ and corresponding eigenvalues $\lambda_j$ with
\[
	\Delta_g e_j = -\lambda_j^2 e_j.
\]
By this convention, $\lambda_j$ is really the eigenvalue of the half-Laplacian $\sqrt{-\Delta_g}$, and is sometimes called the frequency of $e_j$. We will let
\[
	E_\lambda = \operatorname{span} \{ e_j : \lambda_j = \lambda \}
\]
denote the $\lambda$-eigenspace of $\Delta_g$. The dimension of $E_\lambda$ is finite for each $\lambda$. (We refer the reader to \cite{So14} for a thorough treatment of eigenfunction asymptotics on manifolds, and to do Carmo's standard text \cite{dC} for a treatment of Riemannian manifolds.)

Now let $D$ be a region of positive measure on $M$. We are interested in the pairing of $D$ with a discrete set
\[
	\Lambda \subset \bigcup_{\lambda} E_\lambda
\]
which constitutes an orthogonal basis, Riesz basis, or frame for $L^2(D)$. The following result can be viewed as an analog of Fuglede's \cite{Fu} theorem which says that if a subset of ${\Bbb R}^d$ tiles by a lattice, then the dual lattice generates an orthogonal exponential basis for the corresponding Hilbert space. For this, we need the following definition.

\begin{definition} \label{manifoldtiling}
Let $G$ be a subset of the group of isometries of $M$. We say \emph{$D$ tiles $M$ by $G$} if (1) $\phi_1 D \cap \phi_2 D$ has measure zero for all $\phi_1,\phi_2 \in G$ with $\phi_1 \neq \phi_2$, and (2)
\[
    Z \cup \bigcup_{\phi \in G} \phi D = M 
\]
where $Z$ is a set of measure $0$.
\end{definition}

Since $M$ is compact and $D$ has positive measure, the definition implies
\[
    \# G |D| = |M|,
\]
and hence $G$ is finite. Moreover note that a fundamental domain $D$ of a properly discontinuous action by a group $G$ tiles $M$ by $G$.

\begin{theorem} \label{riemannlatticetheorem}
	Let $D \subset M$ be a set of positive measure which tiles $M$ by a subgroup $G$ of the isometries of $M$. Then, there exists an orthogonal basis of eigenfunctions for $L^2(D)$.
\end{theorem}

Theorem \ref{riemannlatticetheorem} relates to the standard Euclidean setting through the flat torus, $\R^d/\Z^d$. Here we choose our eigenfunctions to be the exponentials
\[
    \{ e^{2\pi i m\cdot x} : m \in \Z^d \},
\]
and for $G$ we take any group generated by finitely many rational translations. The proof of Theorem \ref{riemannlatticetheorem} ensures the subset of $G$-periodic exponentials form an orthogonal basis for $L^2(D)$ where $D$ is any fundamental domain of the group action.

\begin{figure}
    \centering
    \includegraphics[scale=3]{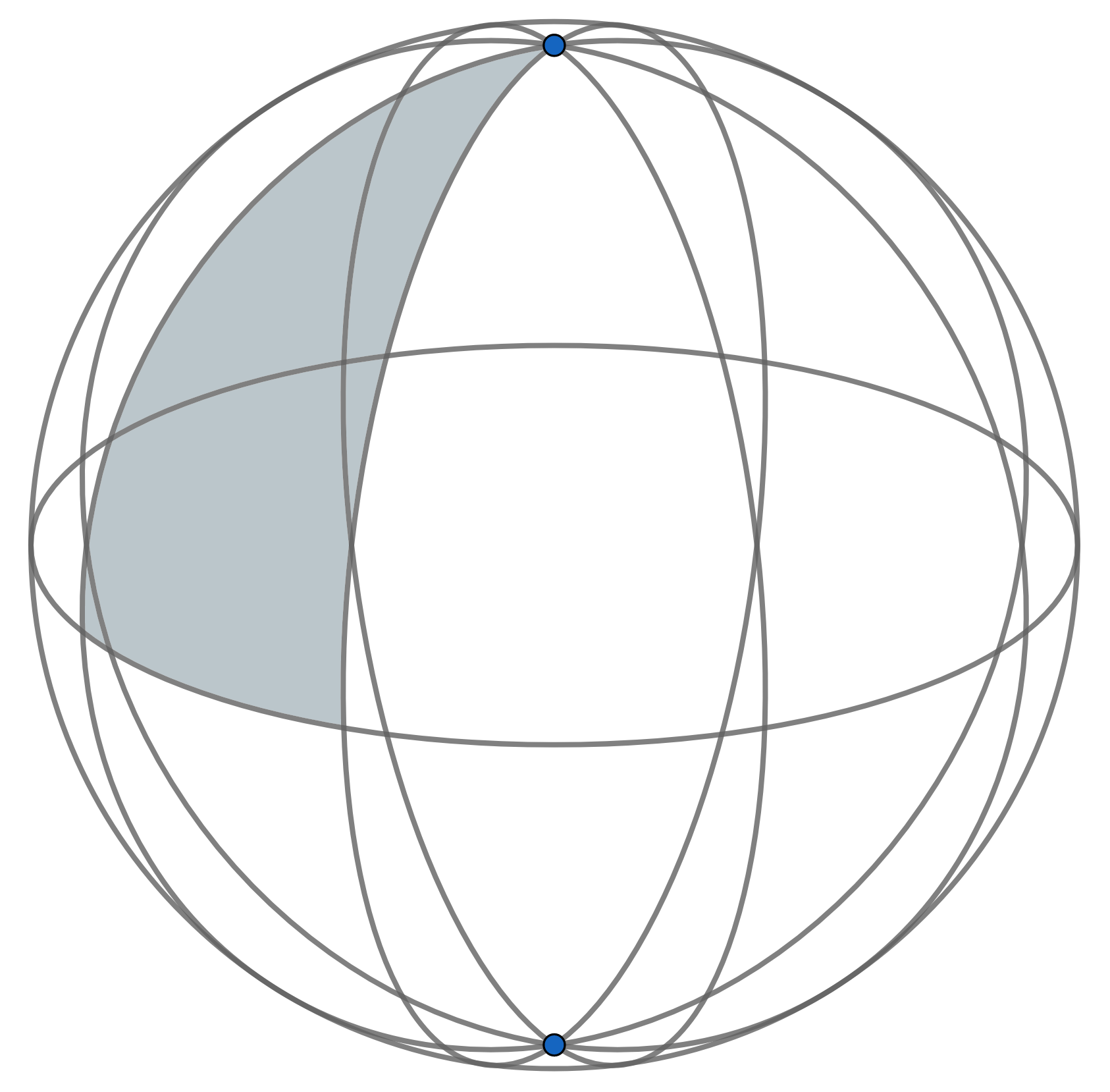}
    \caption{The shaded region tiles $S^2$ by a dihedral group action.}
    \label{spheretilingset}
\end{figure}

The $2$-dimensional sphere $S^2 = \{(x,y,z) : x^2 + y^2 + z^2 = 1\}$ provides a small wealth of examples for which Theorem \ref{riemannlatticetheorem} applies for a noncommutative group action on $M$. For any integer $n \geq 3$, for example, we may consider the action of the dihedral symmetry group $D_n$ on $S^2$. The dihedral group $D_n$ is the group of $2n$ elements generated by a rotation $\sigma$ by $2\pi/n$ about the vertical axis and a rotation by $\pi$ which permutes the north and south poles. Here specifically,
\[
    \sigma = \begin{bmatrix}
        \cos(2\pi/n) & -\sin(2\pi/n) & 0 \\
        \sin(2\pi/n) & \cos(2\pi/n) & 0 \\
        0 & 0 & 1
    \end{bmatrix}
    \qquad \text{ and } \qquad
    \tau = \begin{bmatrix}
        1 & 0 & 0 \\
        0 & -1 & 0 \\
        0 & 0 & -1
    \end{bmatrix}.
\]
Note $\sigma \tau = \tau \sigma^{-1}$, and hence $D_n$ is noncommutative since $\sigma \neq \sigma^{-1}$ for $n \geq 3$. For a set which tiles $S^2$ by the dihedral group action, we take a wedge of the upper hemisphere which makes an angle of $2\pi/n$ at the north pole (figure \ref{spheretilingset}). Similar examples may be found for other symmetry groups of the sphere, e.g. tetrahedral, octahedral, icosahedral, and their respective unoriented versions which include a reflection about the origin.

As remarked above, $G$ need not be abelian, so Theorem \ref{riemannlatticetheorem} is not covered by the generalization of the aforementioned Fuglede's theorem to locally compact abelian groups (see e.g. \cite{BHM}). This result is a small step towards the generalization of the theory of exponential bases and frames to the setting of Riemannian manifolds. We shall study this problem systematically in a sequel.

\medskip

We organize the paper as follows: In Section 2, Theorem \ref{thdecay} and Theorem \ref{thfinite} will be proved. Theorem \ref{maintheorem1} will be proved in Section 3 and we prove our Theorem \ref{maintheorem2} in Section 4. Finally, we establish Theorem \ref{riemannlatticetheorem} in Section 5.

\medskip

\section{Proof of Theorem \ref{thdecay} and \ref{thfinite}}

Throughout the paper, the Fourier transform of a finite Borel measure is defined to be 
$$
\widehat{\mu}(\xi) = \int e^{-2\pi i \xi\cdot x}d\mu(x).
$$
We will also use the notation $X\lesssim Y$ to denote $X\le C Y$ for some constant $C$ that is independent of the variables that defines $X, Y$ and  $X\approx Y$ to denote $X\lesssim Y$ and $X\gtrsim Y$. The measure has a {\it polynomial decay at infinity of order $\gamma$} if the following holds:
\begin{equation}\label{decay}
|\widehat{\mu}(\xi)| \lesssim |\xi|^{-\gamma/2}, \  \forall |\xi|\ge 1.
\end{equation}

\medskip
\noindent {\bf Proof of Theorem \ref{thdecay}}.
We will argue by contradiction. Suppose
$$
\sum_{\lambda\in\Lambda\setminus\{0\}}\frac{1}{|\lambda|^{\gamma}} < \infty.
$$

Take $f(x)=e^{2\pi i \xi\cdot x}$ in \eqref{frame_lower_bound}. It follows that
$$
A \le \sum_{\lambda\in\Lambda}|\widehat{\mu}(\lambda+\xi)|^2
$$
uniformly for all $\xi\in{\mathbb R}^d$. For any $R>1$, and any $|\lambda|>2R$, $|\xi| \le R$, we have  $|\lambda+\xi| > |\lambda|/2$. Therefore 
$$
\sum_{|\lambda|>2R}|\widehat{\mu}(\lambda+\xi)|^2 \lesssim \sum_{|\lambda|>2R}|\lambda+\xi|^{-\gamma} \lesssim \sum_{|\lambda|>2R}|\lambda|^{-\gamma}
$$
As the sum is finite, we can take  $R$ large enough so that 
$$
\sum_{|\lambda|>2R}|\widehat{\mu}(\lambda+\xi)|^2 <\frac{A}{2}.
$$
Then for any $|\xi|\le R$ we have
$$
\frac{A}{2} \le \sum_{|\lambda|\le 2R}|\widehat{\mu}(\lambda+\xi)|^2.
$$

We now integrate this inequality in $\xi$ over the ball $B_R(\vec{0})$, and obtain
$$
\begin{aligned}
 R^d \lesssim& \sum_{|\lambda|\le 2R}\int_{|\xi|\le R}|\widehat{\mu}(\lambda+\xi)|^2 d\xi
 \\
 = & \sum_{|\lambda|\le 2R}\int_{B_R(\lambda)}|\widehat{\mu}(\xi)|^2 d\xi \\
 \le &  \sum_{|\lambda|\le 2R}\int_{B_{3R}(\vec{0})}|\widehat{\mu}(\xi)|^2 d\xi  \  \ \ (\mbox{because $B_R(\lambda)\subset B_{3R}(\vec{0})$})\\
  = & \#(\Lambda \cap B_{2R}(\vec{0}))\cdot \int_{B_{3R}(\vec{0})}|\widehat{\mu}(\xi)|^2 d\xi.
\end{aligned}
$$
Applying the decay condition of $\widehat{\mu}$, it follows that 
\begin{equation}\label{eq_lower_0}
\int_{B_{3R}(\vec{0})}|\widehat{\mu}(\xi)|^2 d\xi \lesssim \int_1^{3R} r^{-\gamma} r^{d-1}dr \lesssim \left\{\begin{array}{ll} R^{d-\gamma} & \hbox{if $\gamma<d$}\\ \log R & \hbox{if $\gamma=d$}\end{array}\right.
\end{equation}
 (Here $\log$ is the natural logaritheorem). This implies that we can find a constant $c$, depending only on $d$, $\gamma$, and $A$, such that for all $R$ large enough,
\begin{equation}\label{lower_bound}
 \#(\Lambda\cap B_{R}(\vec{0})) \ge c\cdot \left\{\begin{array}{ll} R^{\gamma} & \hbox{if $\gamma<d$}\\ \frac{R^d}{\log R} & \hbox{if $\gamma=d$}\end{array}\right.. 
\end{equation}
We finally claim that (\ref{lower_bound}) actually implies that $\sum_{\lambda\in\Lambda\setminus\{0\}}|\lambda|^{-\gamma} = \infty$, from which we obtain our desired contradiction. 

\bigskip

%Indeed, we let $\Gamma = \sum_{\lambda\in\Lambda}\delta_{\lambda}$, where $\delta_{\lambda}$ is the Dirac point mass measure at $\lambda$.
By decomposing the sum into annuli regions and applying the Abel summation formula, we deduce that 
$$
\sum_{|\lambda|>R} \frac{1}{|\lambda|^{\gamma}} \gtrsim \int_R^{\infty} \#(\Lambda\cap B_r(\vec{0})) \frac{1}{r^{\gamma+1}}dr.
$$
Then \eqref{lower_bound} implies
$$
\sum_{|\lambda|>R} \frac{1}{|\lambda|^{\gamma}} \gtrsim \left\{\begin{array}{ll} \int_R^{\infty} \frac1{r}dr & \hbox{if $\gamma<d$}\\ \int_R^{\infty}\frac{1}{r\log r}dr & \hbox{if $\gamma=d$}\end{array}\right.. 
$$
In both cases, the sum diverges. This completes the proof. \qquad$\Box$

\medskip

\noindent {\bf Proof of Theorem \ref{thfinite}}.
Since $E(\Lambda)$ forms a Bessel sequence, 
$$
\sum_{\lambda\in \Lambda, |\lambda|>L} |\widehat{\mu}(\xi+\lambda)|^2\le B.
$$
By integrating both sides in $\xi$ over $B(0,r)$, we obtain
$$
\sum_{\lambda\in \Lambda, |\lambda|>L} \int_{B_r(\lambda)} |\widehat{\mu}(\xi)|^2d\xi =\sum_{\lambda\in \Lambda, |\lambda|>L} \int_{B_r(\vec{0})} |\widehat{\mu}(\xi+\lambda)|^2d\xi \lesssim B r^d.
$$
Invoking the definition of $c$ in Theorem \ref{thfinite}, it follows that 
$$
c \cdot \sum_{\lambda\in \Lambda, |\lambda|>L} \frac{1}{|\lambda|^{\gamma}} \lesssim r^d,
$$
which shows that $\sum_{\lambda\in \Lambda\setminus\{0\}} \frac{1}{|\lambda|^{\gamma}}<\infty$. This completes the proof. \qquad$\Box$

\medskip

\section{Proof of Theorem \ref{maintheorem1}}
Let $K$ be a compact, convex body with smooth boundary $\partial K$ of positive Gaussian curvature. Let $\rho$ be the Minkowski functional associated to $K$ so that $K = \{x\in{\mathbb R}^d: \rho(x)\le 1\}$. The dual norm of $\rho$ is given by 
$$
\rho_{\ast}(\xi): = \sup_{x\in \partial K} x\cdot \xi
$$  
Let $\sigma$ be the surface area measure on $\partial K$. The following Fourier asymptotic formula of $\widehat{\sigma}$ was proved by Herz \cite{Herz}.

\medskip

\begin{theorem}\label{Herz}
Let $K$ be a convex body on ${\mathbb R}^d$ with $\partial K$ smooth and everywhere positive Gaussian curvature. Then 
$$
\widehat{\sigma}(\xi) = C\left(\frac{\xi}{|\xi|}\right)|\xi|^{-\frac{d-1}{2}} \cos\left( 2\pi \left( \rho_{\ast}(\xi)-\frac{d-1}{8}\right)\right) + D_K(\xi)
$$
where $C$ is some positive continuous function and $|D_K(\xi)|\lesssim |\xi|^{-\frac{d+1}{2}}$.  
\end{theorem}
\medskip

\noindent{\bf Proof of Theorem \ref{maintheorem1}.}  Suppose the surface measure $\sigma$ on $\partial K$ is frame-spectral with a spectrum $\Lambda$. We will show that all conditions in Theorem \ref{thdecay} and Theorem \ref{thfinite} are satisfied. Then contradiction follows on the summability of $\sum_{\lambda\in\Lambda\setminus\{0\}} |\lambda|^{-(d-1)}$.

\medskip

By Theorem \ref{Herz}, $\widehat \sigma$ decays polynomially of order $d-1$, so conditions in Theorem \ref{thdecay} are satisfied. To finish the proof, it suffices to show that (\ref{eqthfinite}) holds for $\widehat{\sigma}$.

\medskip

Let $r=r(K)>0$ be a fixed constant that will be determined later. Apply Theorem \ref{Herz} to $\sigma$. Since $|C(\xi/|\xi|)|\approx 1$, by the inequality $|a+b|^2\ge \frac12|a|^2-|b|^2$ we have
$$
\begin{aligned}
|\lambda|^{d-1}\cdot \int_{B_r(\lambda)} |\widehat{\sigma}(\xi)|^2d\xi \gtrsim& |\lambda|^{d-1}\cdot\int_{B_r(\lambda)} |\xi|^{-(d-1)} \cos^2\left( 2\pi \left( \rho_{\ast}(\xi) - \frac{d-1}{8}\right) \right) d\xi \\ & \ \  - |\lambda|^{d-1}\cdot \int_{B_r(\lambda)} |D_K(\xi)|^2d\xi\\
: =& I-II. 
\end{aligned}
$$
We shall show that $I\geq c_K>0$ while $II$ is a small error.

\medskip
Recall that $r$ is a fixed constant, so we may assume $L>2r$. Since we only work with $\xi\in B_r(\lambda)$, $|\lambda|>L$, it follows that $|\xi|\approx |\lambda|$, and therefore
$$
II\lesssim |\lambda|^{d-1} \cdot \int_{B_r(\lambda)}\frac{1}{|\xi|^{d+1}}d\xi \leq \frac{C_r}{|\lambda|}
$$
where $C_r$ is some constant depending only on $r$. This means, for all $\epsilon>0$, we can find $L=L(r)$ large enough so that $ II  <\epsilon$ whenever $|\lambda|>L$.

\medskip
 
Now it remains to show $ I \ge c_K>0$ uniformly in $\lambda$ when $|\lambda|\ge L$ for some $L$. As we just discussed, $|\lambda|\approx |\xi|$, so
\begin{equation}\label{eqI}
I \gtrsim \int_{B_r(\lambda)} \cos^2\left( 2\pi \left( \rho_{\ast}(\xi) - \frac{d-1}{8}\right)\right) d\xi.
\end{equation}

We make the following geometric observation that will be proved at the end of this section.

\medskip

\noindent{\bf Claim:}
	There exists $r=r(K), c=c(K)>0$ such that for any $\lambda\in\R^d$ such that $|\lambda|$ and $ \rho_*(\lambda)>100$, for any $\rho_*(\lambda)-1<t< \rho_*(\lambda)+1$, the cap
	$$B_r(\lambda)\cap\{\xi: \rho_*(\xi)=t\}$$
	has diameter $\geq c>0$.
	
\medskip

Notice that $\rho_*$ is homogeneous of degree $1$, so the spherical cap $B_r(\lambda)\cap\{\xi: \rho_*(\xi)=t\}$ in the claim has $(d-1)$-dimensional Hausdorff measure
$$t^{d-1}\cdot {\mathcal H}^{d-1}(\text{a $ct^{-1}$-cap on }\rho_*=1)\gtrsim_K 1.$$
Here, $\gtrsim_K$ means that the constant depends only on the convex body $K$. We now assume that the claim holds.  Since $\rho_*$ is differentiable (see e.g., Corollary 25.1.3, \cite{Roc}), it follows that $\nabla \rho_*$ is homogeneous of degree $0$ and $|\nabla \rho_*|\approx 1$ on $\R^d\backslash\{0\}$. Then by the coarea formula, which we may think it is the polar coordinates with respect to the $\rho_*$-metric, together with our  claim, we have
$$
\begin{aligned}
I\gtrsim& \int_{\rho_{\ast}(\lambda)-1}^{\rho_{\ast}(\lambda)+1}\left(\int_{B_r(\lambda)\cap\{\xi: \rho_*(\xi)=t\}} \frac{1}{|\nabla \rho_*(\xi)|}\,d{\mathcal H}^{d-1}(\xi)\right)\cos^2\left( 2\pi \left(t - \frac{d-1}{8}\right)\right) dt \\
\gtrsim_K &\int_{\rho_{\ast}(\lambda)-1}^{\rho_{\ast}(\lambda)+1} \cos^2\left( 2\pi \left(t - \frac{d-1}{8}\right)\right) dt.\\
\end{aligned}
$$

Notice that the interval $\left[\rho_{\ast}(\lambda)-1,\rho_{\ast}(\lambda)+1\right]$ always contains an interval $[k-1/2,k+1/2]$ where $k$ is a positive integer. By the integral periodicity of $\cos2\pi x$,
$$
I\gtrsim_K \int_{-1/2}^{1/2} \cos^2\left( 2\pi \left(t - \frac{d-1}{8}\right)\right) dt =: c_K>0
$$
and $c_K$ is a constant independent of $\lambda$.  To finish the proof, we need to show \eqref{eqthfinite}. First we fix $r=r(K)$ in the claim above, then when $L=L(r, K)$ is large enough, $I\geq c_K>0$ while $II<c_K/2$ for any $|\lambda|>L$, as desired. 

\medskip

It remains to justify our claim.

\medskip

\noindent {\it Proof of the claim:}  Denote $|\cdot|$ as the Euclidean norm. By the convexity of $K$, there exists $C_K\ge 1$ such that 
$$
C_{K}^{-1} |x| \le \rho_{\ast}(x) \le C_K |x|, \ \forall x\in{\mathbb R}^d.
$$
For any $t\in(\rho_*(\lambda)-1, \rho_*(\lambda)+1)$, we define $\lambda_t = \frac{t}{\rho_{\ast}(\lambda)} \lambda$. Then $\rho_{\ast}(\lambda_t) = t$ by the homogenity of $\rho_{\ast}$, and
$$
|\lambda-\lambda_t|\leq C_K\cdot \rho_*(\lambda-\lambda_t) =  C_K\cdot \left|1-\frac{t}{\rho_{\ast}(\lambda)}\right|\cdot\rho_{\ast}(\lambda) = C_K\cdot |\rho_{\ast}(\lambda)-t|\leq C_K. 
$$
This implies, if we choose $r=2 C_K$, then for any $\rho_*(\lambda)-1\leq t\leq \rho_*(\lambda)+1$,
$$B_{r/2}(\lambda_t)\subset B_r(\lambda).$$

Given any $\omega$ such that $\rho_*(\omega)=1$ and any small enough $\delta>0$, for example $\delta<\frac{1}{100}\cdot \mbox{diam}(\{\xi: \rho_{\ast}(\xi) = 1\})$, the diameter of $B(\omega, \delta)\cap \{\rho_*=1\}$ is $\geq c_K \delta>0$ for some constant $c_K>0$. Hence,  by the homogeneity of $\rho_*$, for any $\rho_*(\lambda)-1<t< \rho_*(\lambda)+1$, the diameter of $B_{r}(\lambda)\cap \{\xi: \rho_*(\xi)=t\}$ is at least
$$\mbox{diam} \left(B_{r/2}(\lambda_t))\cap \{\xi: \rho_*(\xi)=t\}\right)\gtrsim|\lambda_t| \cdot c_K |\lambda_t|^{-1}=c_K>0,$$
as desired.\qquad$\Box$

\medskip

\section{Fourier frame on surfaces without curvature}

We will prove Theorem \ref{maintheorem2} in this section. In particular, we prove the following more general theorem.

\begin{theorem}\label{theorem_flat_surface}
	Let ${\mathcal F}$ be a collection of finitely many bounded Borel subsets of $\R^d$ such that 
	\begin{itemize}
		\item each $F\in{\mathcal F}$ is contained in a $k_F$-dimensional affine subspace $V_F\subset\R^d$, ${\mathcal H}^{k_F}(F)>0$, and
		\item for any $F, F'\in {\mathcal F}$, $k_F> k_{F'}$, we have $V_{F'}+\tau\not\subset V_F$ for any $\tau\in\R^d$.
	\end{itemize}

	Then the measure
	$$\sum_{F\in{\mathcal F}}{\mathcal H}^{k_F}\big|_{F} $$
	admits a Fourier frame.
\end{theorem}

\medskip

To see why Theorem \ref{maintheorem2} follows from Theorem \ref{theorem_flat_surface}, we note that if $\sigma$ is the surface measure of a polytope, then 
$$
\sigma = \sum_{F\in{\mathcal F}}{\mathcal H}^{d-1}\big|_{F} 
$$
where $F$ are all $(d-1)$-dimensional facet of the polytope. Hence, the first condition is satisfied and the second condition is vacuously true. Therefore, Theorem \ref{theorem_flat_surface} applies.

\begin{remark}
{\textrm The second condition in Theorem \ref{theorem_flat_surface} cannot be removed. One can consider  a union of a line $L$ and a square $Q$ in the same plane, so that the second condition is violated. By regarding the plane as ${\mathbb R}^2$, the measure ${\mathcal H}^{1}\big|_{L}+{\mathcal H}^{2}\big|_{Q}$ will then be a non-trivial sum of a singular measure ${\mathcal H}^{1}\big|_{L}$ and an absolutely continuous measure ${\mathcal H}^{2}\big|_{Q}$. By the law of pure type of frame-spectral measure \cite{HLL}, this measure does not admit any Fourier frame.}
\end{remark}

\medskip

We now start to prove Theorem \ref{theorem_flat_surface}. Since we can always put a set $F$ into a cube and a Fourier frame on $F$ will induce naturally a Fourier frame for its subset, we may without loss of generality assume each $F\in{\mathcal F}$ is a $k_F$-dimensional unit cube in $V_F$, and $V_F\neq V_{F'}$ for any $F, F'\in\mathcal{F}$. Moreover, we assume any two such unit cubes contained in parallel spaces are translates of each other.

\medskip

We now divide ${\mathcal F}$ into equivalent classes 
$$
{\mathcal F}_1,\dots, {\mathcal F}_m,$$
 where $F\sim F'$ if and only if there exists $\tau\in\R^d$ such that $F+\tau=F'$. For each class ${\mathcal F}_j$ , there exists a $k_j$-dimensional subspace $V_j$, a $k_j$-dimensional unit cube $Q_j\subset V_j\subset\R^d$, translations $\tau_j^1,\dots,\tau_j^{\#({\mathcal F}_j)}\in\R^d $ such that
$$
{\mathcal F}_j= Q_j +\{\tau_j^1,\dots,\tau_j^{\#({\mathcal F}_j)}\}.
$$
We note that from our equivalent class definition and also the second assumption stated in Theorem \ref{theorem_flat_surface}, none of the  $V_j$ is contained in the other $V_j'$. 
Denote $\Z^{k_j}\subset V_j\approx\R^{k_j}$ as a natural spectrum of $Q_j$. Also denote by  $P_{V_j}:\R^d\rightarrow V_j$  the orthogonal projection onto the subspace $V_j$. Recall also that a discrete set $\Lambda$ is called {\it $\delta$-separated} if 
$$
|\lambda-\lambda_t|\ge \delta>0,
$$
for all $\lambda\ne\lambda_t\in\Lambda$.
We will need the following lemmas. The proof of these lemma will be given in the end of the section. 

\medskip

\begin{lemma}\label{non-parallel-faces}
	For any $N>0$ and $ \delta>0$, there exist discrete sets $\Gamma_j\subset\R^d$, $j=1,\dots,m$ such that  
\begin{enumerate}
    \item $\#\left(\Gamma_j \cap P_{V_j}^{-1}(z)\right)=N$ for any $z\in {\mathbb Z}^{k_j}$, $j=1,\dots, m$, and
	\item  for any pair $j, \ell$	, $j\neq \ell$, the map $P_{V_{\ell}}: \Gamma_j\rightarrow V_{\ell}$ is injective and $P_{V_{\ell}}(\Gamma_j)$ is $\delta$-separated in $V_{\ell}$.
\end{enumerate}
\end{lemma}

\medskip

\begin{lemma}\label{parallel-faces}
	Under notation above, there exist $\epsilon>0$ such that for all $j=1,...,m$,  we can find finite sets 
	$$
	A_j=\{\alpha_j^1,\dots,\alpha_j^{\#({\mathcal F}_j)}\}\subset V_j^\perp\subset\R^d,
	$$ such that 
	$$ \epsilon\le \|{\bf M}_j\| \ $$%\mbox{and} \ \|{\bf M}_j\|^{-1}<\epsilon^{-1}, $$
	where ${\bf M}_j$ is the $\#({\mathcal F}_j)\times \#({\mathcal F}_j)$ matrix $\left(e^{-2\pi i \tau_j^{\ell}\cdot \alpha^s_j}\right)$, $\ell, s=1,\dots, \#({\mathcal F}_j)$ and $\|\cdot\|$ denotes the matrix lower bound $\|M\| = \inf_{x \neq 0} |Mx|/|x|$.
\end{lemma}

\medskip

\begin{lemma}\label{Young-book}
	Let $\Gamma\subset\R^d$ be discrete and $\delta$-separated. Then 
	$$E(\Gamma):=\{e^{2\pi i \gamma\cdot x}: \gamma\in\Gamma\}
	$$ is a Bessel sequence  for $L^2([0,1]^d)$, where the upper Bessel bound depends only  on $\delta$.
\end{lemma}

\noindent{\bf Proof of Theorem \ref{theorem_flat_surface}.} Let $\Gamma_j$ be the discrete set in Lemma \ref{non-parallel-faces} with $N$ to be determined and let $A_j$ be the finite sets determined in Lemma \ref{parallel-faces}. We then define
 $$
 \Lambda_j = \Gamma_j+A_j.
 $$
  We shall show that when $N$ is large enough,
$$
\left\{ e^{2\pi i \lambda\cdot x}: \lambda\in\Lambda:=\bigcup_{i=1}^m \Lambda_i\right\}$$
is a Fourier frame of $$\sum_{F\in{\mathcal F}}{\mathcal H}^{k_F}\big|_{F}=\sum_{j=1}^m \sum_{\ell=1}^{\#({\mathcal F}_j)}{\mathcal H}^{k_j}\big|_{Q_j+\tau_j^{\ell}}.
$$
Let $f_{j\ell}\in L^2(Q_j+\tau_j^{\ell})$, $j=1,\dots, m$, $\ell=1,\dots, \#({\mathcal F}_j)$. Consider
\begin{equation}\label{l2-sum}
\sum_{\lambda\in\Lambda} \left|\sum_{j=1}^m \sum_{\ell=1}^{\#({\mathcal F}_j)}\int_{Q_j+\tau_j^{\ell}}e^{-2\pi i \lambda\cdot x} f_{j\ell}(x)\,d{\mathcal H}^k(x)\right|^2 . 
\end{equation}
The upper frame bound follows easily by Lemma \ref{non-parallel-faces} and Lemma \ref{Young-book}. So it suffices to show the lower frame bound, namely
$$
\sum_{\lambda\in\Lambda} \left|\sum_{j=1}^m \sum_{\ell=1}^{\#({\mathcal F}_j)}\int_{Q_j+\tau_j^{\ell}}e^{-2\pi i \lambda \cdot x} f_{j{\ell}}(x)\,d{\mathcal H}^k(x)\right|^2\ge C_{N, m, \delta, \epsilon} \sum_{j=1}^m\sum_{\ell=1}^{\#({\mathcal F}_j)} ||f_{j\ell}||_{L^2(Q_j+\tau_j^{\ell})}^2. 
$$
for some constant $C_{N,m,\delta,\epsilon}$ that depends only on $N,m,\delta,\epsilon$ but not $f$. 
We now decompose the sum $\sum_{\lambda\in\Lambda}$ in \eqref{l2-sum} into $\sum_{i=1}^m\sum_{\lambda\in\Lambda_i}$. We first estimate the sum over $\lambda\in\Lambda_1$. Then by the inequality
$$(a_1+\cdots+a_m)^2\geq \frac{1}{m} a_1^2-a_2^2-\cdots-a_{m}^2,$$
the sum in \eqref{l2-sum}, with $\lambda$ only in $\Lambda_1$,  is bounded from below by
$$
\begin{aligned}
&\sum_{\lambda \in\Lambda_1 }\left(\frac{1}{m}\left|\sum_{\ell=1}^{\#({\mathcal F}_1)}\int_{Q_1+\tau_1^{\ell}}e^{-2\pi i x\cdot\lambda} f_{1\ell}(x)\,d{\mathcal H}^k(x)\right|^2-\sum_{j=2}^m\left|\sum_{\ell=1}^{\#({\mathcal F}_j)}\int_{Q_j+\tau_j^{\ell}}e^{-2\pi i \lambda\cdot x} f_{j\ell}(x)\,d{\mathcal H}^k(x)\right|^2\right)\\
=& \frac{1}{m}\sum_{\lambda \in\Lambda_1 }\left|\sum_{\ell=1}^{\#({\mathcal F}_1)}\int_{Q_1+\tau_1^{\ell}}e^{-2\pi i \lambda\cdot x} f_{1\ell}(x)\,d{\mathcal H}^k(x)\right|^2-\sum_{\lambda \in\Lambda_1}\sum_{j=2}^m\left|\sum_{\ell=1}^{\#({\mathcal F}_j)}\int_{Q_j+\tau_j^{\ell}}e^{-2\pi i \lambda \cdot x} f_{j\ell}(x)\,d{\mathcal H}^k(x)\right|^2\\
:=& I-II.
\end{aligned}
$$
Since $\Lambda_1=\Gamma_1+A_1$, $A_1\subset V_1^\perp$, we can rewrite $I$ as
\begin{equation*}
\begin{aligned}I=&\frac{1}{m}\sum_{\gamma\in\Gamma_1}\sum_{\alpha\in A_1}\left|\sum_{\ell=1}^{\#({\mathcal F}_1)} e^{-2\pi i \tau_1^{\ell}\cdot\alpha}\left(e^{-2\pi i \tau_1^{\ell}\cdot \gamma}\int_{Q_1} e^{-2\pi i (\gamma+\alpha)\cdot x} f_{1\ell}(x+\tau_1^{\ell})\,d{\mathcal H}^k(x)\right)\right|^2\\
=&\frac{1}{m} \sum_{\gamma\in\Gamma_1}\sum_{\alpha\in A_1}\left|\sum_{\ell=1}^{\#({\mathcal F}_1)} e^{-2\pi i \tau_1^{\ell}\cdot\alpha}\left(e^{-2\pi i \tau_1^{\ell}\cdot \gamma}\int_{Q_1} e^{-2\pi i \gamma \cdot x} f_{1\ell}(x+\tau_1^{\ell})	\,d{\mathcal H}^k(x)\right)\right|^2\\
=& \frac{1}{m} \sum_{\gamma\in\Gamma_1}\sum_{\alpha\in A_1}\left|\sum_{\ell=1}^{\#({\mathcal F}_1)} e^{-2\pi i \tau_1^{\ell}\cdot\alpha} v_{\ell}\right|^2
\end{aligned}
\end{equation*}
where $v_{\ell} = e^{-2\pi i \tau_1^{\ell}\cdot \gamma}\int_{Q_1} e^{-2\pi i \gamma \cdot x} f_{1\ell}(x+\tau_1^{\ell})	\,d{\mathcal H}^k(x)$. Letting ${\bf v} = (v_1,...,v_{\#{\mathcal F}_1})^{\mathtt T}$ , we note that the inner sum is equal to $\|M_1{\bf v}\|^2$, where
 ${\mathbf M}= \left(e^{-2\pi i \tau_1^l\cdot \alpha^s_1}\right)$ as defined in  Lemma \ref{parallel-faces}. Using Lemma \ref{parallel-faces}, the matrix lower bound $\|{\mathbf M}\| \ge \epsilon$. % $\left|\left|\left(e^{-2\pi i \tau_1^l\cdot \alpha^s_1}\right)^{-1}\right|\right|>\epsilon$.
 It follows that
$$
I\geq \frac{\epsilon^2}{m} \sum_{\gamma\in\Gamma_1}\sum_{\ell=1}^{\#({\mathcal F}_1)}\left|\int_{Q_1} e^{-2\pi i \lambda\cdot x} f_{1\ell}(x+\tau_1^{\ell})	\,d{\mathcal H}^k(x)\right|^2,
$$
which by (1) in Lemma \ref{non-parallel-faces} equals
$$\frac{\epsilon^2 N}{m}\sum_{\ell=1}^{\#({\mathcal F}_1)} ||f_{1\ell}||_{L^2(Q_1+\tau_1^{\ell})}^2.$$

\medskip

Now we estimate $II$.  we can rewrite $II$ as
\begin{equation*}
\begin{aligned}
II=&\sum_{j=2}^m\sum_{\gamma\in\Gamma_1}\sum_{\alpha\in A_1}\left|\sum_{\ell=1}^{\#({\mathcal F}_j)} e^{-2\pi i \tau_j^{\ell}\cdot (\gamma+\alpha)}\left(\int_{Q_j} e^{-2\pi i (\gamma+\alpha)\cdot x} f_{j{\ell}}(x+\tau_j^{\ell})\,d{\mathcal H}^k(x)\right)\right|^2\\
\le & \sum_{j=2}^m\sum_{\gamma\in\Gamma_1}\sum_{\alpha\in A_1}(\#{\mathcal F}_j)\sum_{\ell=1}^{\#({\mathcal F}_j)}\left|\int_{Q_j} e^{-2\pi i (\gamma+\alpha)\cdot x} f_{j{\ell}}(x+\tau_j^{\ell})\,d{\mathcal H}^k(x)\right|^2\\ 
\le& \left(\max \#({\mathcal F}_j)\right) \sum_{j=2}^m \sum_{\alpha\in A_1}\sum_{\ell=1}^{\#({\mathcal F}_j)}\sum_{\gamma\in\Gamma_1}\left|\int_{Q_j} e^{-2\pi i (P_{V_j}(\gamma))\cdot x} \left(e^{-2\pi i \alpha\cdot x}f_{j{\ell}}(x+\tau_j^{\ell})\right)\,d{\mathcal H}^k(x)\right|^2\\
%=& \sum_{j=2}^m\sum_{\gamma\in\Gamma_j}\sum_{\alpha\in A_j}\left|\sum_{\ell=1}^{\#({\mathcal F}_j)} e^{-2\pi i \tau_j^{\ell}\cdot\alpha}\left(e^{-2\pi i \tau_j^{\ell}\cdot \gamma}\int_{Q_j} e^{-2\pi i \gamma\cdot x} f_{j{\ell}}(x+\tau_j^{\ell})\,d{\mathcal H}^k(x)\right)\right|^2
\end{aligned}
\end{equation*}
 where we used a Cauchy-Schwarz inequality in the inner summation. By (2) in Lemma \ref{non-parallel-faces}, $P_{V_j}(\Gamma_1)$ is $\delta$-separated. Therefore, Lemma \ref{Young-book} tells us that 
 $$
 \begin{aligned}
 II\le& \left(\max \#({\mathcal F}_j)\right) \cdot  C_{\delta}\cdot \sum_{j=2}^m \sum_{\alpha\in A_1}\sum_{\ell=1}^{\#({\mathcal F}_j)}\int_{Q_j} \left|e^{-2\pi i \alpha\cdot x}f_{j{\ell}}(x+\tau_j^{\ell})\right|^2\,d{\mathcal H}^k(x)\\
  = &\left(\max \#({\mathcal F}_j)\right) \cdot  C_{\delta}\cdot \sum_{j=2}^m\sum_{\alpha\in A_1}\sum_{\ell=1}^{\#({\mathcal F}_j)} ||f_{j\ell}||_{L^2(Q_j+\tau_j^{\ell})}^2\\
  \le& M^2 \cdot   C_{\delta}\cdot\sum_{j=2}^m\sum_{\ell=1}^{\#({\mathcal F}_j)} ||f_{j\ell}||_{L^2(Q_j+\tau_j^{\ell})}^2\\
 \end{aligned}
 $$
%  Lemma \ref{parallel-faces} implies that the matrix norm $\left|\left|{\bf M}_j\right|\right|<\epsilon^{-1}$. With a similar matrix norm argument used for estimating $I$, it follows that
%$$II\leq \epsilon^{-2} \sum_{j=2}^m \sum_{\gamma\in\Gamma_j}\sum_{\ell=1}^{\#({\mathcal E}_j)}\left|\int_{Q_j} e^{-2\pi i \lambda\cdot x} f_{j\ell}(x+\tau_j^{\ell})	\,d{\mathcal H}^k(x)\right|^2,$$
%which is by (2) in Lemma \ref{non-parallel-faces} and Lemma \ref{Young-book}
%$$
%\leq C_{\delta, \epsilon} \sum_{j=2}^m\sum_{\ell=1}^{\#({\mathcal F}_j)} ||f_{j\ell}||_{L^2(Q_j+\tau_j^{\ell})}^2.
%$$
where $M =\max \#({\mathcal F}_j)$. Hence
$$\eqref{l2-sum}\geq I-II\geq \frac{\epsilon N}{m}\sum_{\ell=1}^{\#({\mathcal F}_1)} ||f_{1{\ell}}||_{L^2(Q_1+\tau_1^{\ell})}^2 - C_{\delta}M^2 \sum_{j=2}^m\sum_{\ell=1}^{\#({\mathcal F}_j)} ||f_{j{\ell}}||_{L^2(Q_j+\tau_j^{\ell})}^2. $$
Similarly we have for any $i=1,\cdots,m$,
$$
\begin{aligned}
&\sum_{\lambda\in\Lambda_i} \left|\sum_{j=1}^m \sum_{\ell=1}^{\#({\mathcal F}_j)}\int_{Q_j+\tau_j^{\ell}}e^{-2\pi i \lambda\cdot x} f_{j\ell}(x)\,d{\mathcal H}^k(x)\right|^2 \\\geq &\frac{\epsilon N}{m}\sum_{\ell=1}^{\#({\mathcal F}_i)} ||f_{i\ell}||_{L^2(Q_i+\tau_i^{\ell})}^2 - C_{\delta} M^2 \sum_{j\neq i}\sum_{\ell=1}^{\#({\mathcal F}_j)} ||f_{j\ell}||_{L^2(Q_j+\tau_j^{\ell})}^2.
\end{aligned} 
$$
Taking the sum in $i$, we have the following lower bound of \eqref{l2-sum}:
$$
\begin{aligned}
&\sum_{\lambda\in\Lambda} \left|\sum_{j=1}^m \sum_{\ell=1}^{\#({\mathcal F}_j)}\int_{Q_j+\tau_j^{\ell}}e^{-2\pi i \lambda\cdot x} f_{j\ell}(x)\,d{\mathcal H}^k(x)\right|^2\\\geq &\left(\frac{\epsilon N}{m}-(m-1)C_{\delta} M^2\right) \sum_{j=1}^m\sum_{\ell=1}^{\#({\mathcal F}_j)} ||f_{j\ell}||_{L^2(Q_j+\tau_j^{\ell})}^2.
\end{aligned} 
$$
Hence, if $N$ is large enough, there will be a positive lower frame bound.
\qquad$\Box$
\medskip

\noindent{\bf Proof of Lemma \ref{non-parallel-faces}.} We just prove the case when $j=1$, others follow from a similar argument. We first claim that there exists a direction $\omega\in S^{d-1}$ such that $\omega\in V_1^\perp$, $\omega\notin V_j^\perp$ for all $j\geq 2$. Suppose the claim is false. Then 
	$$V_1^\perp\subset \bigcup_{j=2}^m V_j^\perp. $$
	In particular, this means that 
	$$
	V_{1}^{\perp} =  \bigcup_{j=2}^m V_j^\perp \cap V_1^{\perp}.
	$$
	Then there is necessarily some $j_0\ge2$ such that
	% $V_{j_0}^\perp \cap V_1^{\perp}$, as a subspace of $V_{1}^{\perp}$, has positive Lebesgue measure, which implies
	$V_1^\perp \subset V_{j_0}^\perp$. But this is not possible since we know none of the $V_j$ is contained in $V_1$. Therefore, the claim holds. 
	
	\medskip

The claim implies that $P_{V_j}(t\omega)$, $t\in\R$ is a straight line inside $V_j$.
We construct $\Gamma_1$ inductively as follows. We enumerate $\Z^{k_1} \subset V_1$ by $\{z_1,z_2,\ldots \}$ and set
$$
    Z_1 = \{z_1+t^1_1\omega,...,z_1+t^1_N\omega\},
$$
where we have chosen $t^1_s$ for $s = 1,\ldots, N$ so that the projections of $Z_1$ onto each of the other subspaces $V_2,\ldots,V_m$ are $\delta$-separated. For the inductive step, set
$$
    Z_k = Z_{k-1} \cup \{z_{k}+t^{k}_1\omega,...,z_{k}+t^k_N\omega\}
$$
where again the $t^1_s$ for $s = 1,\ldots, N$ are chosen so that the projections of $Z_k$ onto each of $V_2,\ldots, V_m$ remain $\delta$-separated. This is achievable since at each stage $Z_{k-1}$ is finite and $P_{V_j}(t\omega)$ for $t\in\R$ is a straight line inside $V_j$ for $j \geq 2$. We set
$$
    \Gamma_1 = \bigcup_{k = 1}^\infty Z_k
$$
and see that it satisfies parts (1) and (2) of the lemma by construction. The complete lemma follows after constructing $\Gamma_2,\ldots, \Gamma_m$ similarly.
\qquad$\Box$
\medskip

\noindent{\bf Proof of Lemma \ref{parallel-faces}}. We first fix $j\in\{1,...,m\}$. By our reduction, $F, F'\in {\mathcal F}_j$ lie in the same $k_j$-dimensional affine subspace if and only if $F=F'$, so $P_{V_j^\perp}(\tau_j^{\ell})\neq P_{V_j^\perp}(\tau_j^{\ell'})$ for any $j$ and any $\ell \neq \ell'$. Then one can choose $\alpha_j^0\in V_j^\perp\backslash\{0\}$, $|\alpha_j^0|>0$ small enough, such that
	$$e^{-2\pi i \tau_j^{\ell}\cdot \alpha_j^0}\neq e^{-2\pi i \tau_j^{\ell'}\cdot \alpha_j^0},\ \forall\ \ell\neq \ell'.$$
	Let $\alpha_j^s = s\cdot \alpha_j^0$, $s=1,\dots,\#({\mathcal F}_j)$. Then ${\bf M}_j=\left(\left(e^{-2\pi i \tau_j^{\ell}\cdot \alpha^0_j}\right)^s\right)$, $\ell, s=1,\dots, \#({\mathcal F}_j)$ is a Vandermonde matrix whose determinant is not zero. Therefore, $\|{\bf M}_j\|>0$ and we finish the proof by taking $\epsilon = \min\{\|{\bf M}_j\|: j=1,...,m\}$.
\qquad$\Box$

\medskip

\noindent{\bf Proof of Lemma \ref{Young-book}.} We believe this lemma is well-known, but we would like to put it here for the sake of self-containment.  This proof here is based on \cite[Lemma 1]{GR}. Given a continuous function $F$, we define $F^{\#}(x) = \sup_{|y-x|\le \delta} |F(y)|$. We now take $\varphi$ to be a Schwartz space function such that $\varphi\equiv1$ on $[0,1]^d$. Then for any $f\in L^2([0,1]^d)$, $f = f\varphi$ on the unit cube and therefore we have 
$$
\widehat{f} = \widehat{f}\ast\widehat{\varphi}
$$
Writing $F = \widehat{f}$, we can deduce easily that $F^{\#}(x) \le \left(|F|\ast \widehat{\varphi}^{\#}\right)(x)$. We therefore obtain from Young's inequality that 
\begin{equation}\label{eq3}
\|F^{\#}\|_2 \le \|\widehat{\varphi}^{\#}\|_1 \|F\|_2
\end{equation}
with $\|\widehat{\varphi}^{\#}\|_1<\infty$ since $\widehat{\varphi}$ is also in the Schwartz space.  Therefore, 
$$
\begin{aligned}
\sum_{\gamma\in\Gamma} \left|\int_{[0,1]^d} f(x)e^{-2\pi i \gamma\cdot x}dx\right|^2=& \sum_{\gamma\in\Gamma} |F(\gamma)|^2\\
 = &\sum_{\gamma\in\Gamma} \frac{1}{C\delta^d}\int_{|x-\gamma|\le \delta} |F(\gamma)|^2dx \\
\le& \sum_{\gamma\in\Gamma} \frac{1}{C\delta^d}\int_{|x-\gamma|\le \delta} |F^{\#}(x)|^2dx\\
\le &\frac{1}{C\delta^d} \int |F^{\#}(x)|^2dx   \ \ \ \  (\mbox{since $\Gamma$ is $\delta$-separated})\\
\le &\frac{\|\widehat{\varphi}^{\#}\|^2_1 }{C\delta^d} \int |F(x)|^2dx    \ \ \ \  (\mbox{by (\ref{eq3})})\\
=  &\frac{\|\widehat{\varphi}^{\#}\|^2_1 }{C\delta^d} \int |f(x)|^2dx.  
\end{aligned}
$$
We see that the constant depends only the $\delta$, but not on $\Lambda$. This completes the proof.

\medskip

\medskip

\section{Riemannian manifolds: Proof of Theorem \ref{riemannlatticetheorem}} 

\vskip.125in

\noindent{\bf Proof of Theorem \ref{riemannlatticetheorem}.} Let $G$ be our subgroup of isometries as in Theorem \ref{riemannlatticetheorem}. Since $D$ tiles $M$ by $G$, we are done provided we construct an orthogonal basis for $L^2(M / G)$, the space of $G$-periodic functions in $L^2(M)$. We define a projection operator $P$ onto $G$-periodic functions by
	\[
		P f = \frac{1}{\#G} \sum_{\phi \in G} f \circ \phi.
	\]
	Note $P E_\lambda \subset E_\lambda$, $P^2 = P$, and $P$ is self-adjoint. We select our orthogonal basis $e_j$ to diagonalize $P|_{E_\lambda}$ for each $\lambda$. Since $P^2 = P$, these $e_j$'s fall into exactly one of two categories, (i) $Pe_j = e_j$, or (ii) $Pe_j = 0$.
	We take $\Lambda$ to consist of those basis elements satisfying (i). $\Lambda$ inherits orthogonality immediately. Moreover if $f$ is $G$-periodic, then
	\[
		f = \sum_j \langle f, e_j \rangle e_j = \sum_j \langle P f, e_j \rangle e_j = \sum_j \langle f, Pe_j \rangle e_j = \sum_{e_j \in \Lambda} \langle f, e_j \rangle e_j,
	\]
	and hence $\Lambda$ spans $L^2(M/G)$.

\end{document}